\documentclass[11pt,reqno]{article}

\usepackage[usenames]{color}
\usepackage{amssymb}
\usepackage{amsmath}
\usepackage{amsthm}
\usepackage{amsfonts}
\usepackage{amscd}
\usepackage{graphicx}

\usepackage[colorlinks=true,
linkcolor=webgreen,
filecolor=webbrown,
citecolor=webgreen]{hyperref}

\definecolor{webgreen}{rgb}{0,.5,0}
\definecolor{webbrown}{rgb}{.6,0,0}
\usepackage{color}
\usepackage{fullpage}
\usepackage{float}

\usepackage{cases}

\usepackage{graphics,amsmath,amssymb}
\usepackage{amsthm}
\usepackage{amsfonts}
\usepackage{latexsym}

\begin{document}

\theoremstyle{plain}
\newtheorem*{theorem}{Theorem}
\newtheorem{remark}{Remark}

\newtheorem*{proposition}{Proposition}

\begin{center}
\vskip 1cm{\Large\bf A note on the hyper-sums of powers of integers,  \\
\vskip .0in hyperharmonic polynomials and $r$-Stirling \\
\vskip .07in numbers of the first kind}
\vskip .2in \large Jos\'{e} Luis Cereceda \\
{\normalsize Collado Villalba, 28400 (Madrid), Spain} \\
\href{mailto:jl.cereceda@movistar.es}{\normalsize{\tt jl.cereceda@movistar.es}}
\end{center}

\begin{abstract}
Recently, Karg{\i}n {\it et al.\/} (arXiv:2008.00284 [math.NT]) obtained (among many other things) the following formula for the hyper-sums of powers of integers $S_k^{(m)}(n)$
\begin{equation*}
S_k^{(m)}(n) = \frac{1}{m!} \sum_{i=0}^{m} (-1)^i \genfrac{[}{]}{0pt}{}{m+n+1}{i+n+1}_{n+1} S_{k+i}(n),
\end{equation*}
where $S_k^{(0)}(n) \equiv S_k(n)$ is the ordinary power sum $1^k + 2^k + \cdots + n^k$. In this note we point out that a formula equivalent to the preceding one was already established in a different form, namely, a form in which $\genfrac{[}{]}{0pt}{}{m+n+1}{i+n+1}_{n+1}$ is given explicitly as a polynomial in $n$ of degree $m-i$. We find out the connection between this polynomial and the so-called $r$-Stirling polynomials of the first kind. Furthermore, we determine the hyperharmonic polynomials and their successive derivatives in terms of the $r$-Stirling polynomials of the first kind, and show the relationship between the (exponential) complete Bell polynomials and the $r$-Stirling numbers of the first kind. Finally, we derive some identities involving the Bernoulli numbers and polynomials, the $r$-Stirling numbers of the first kind, the Stirling numbers of both kinds, and the harmonic numbers.
\end{abstract}

\section{Introduction}

For integers $k,n \geq 0$ and $m \geq 1$, the hyper-sums of powers of integers $S_k^{(m)}(n)$ are defined recursively by
\begin{equation*}
S_k^{(m)}(n) = \sum_{j=1}^{n} S_k^{(m-1)}(j),
\end{equation*}
with initial condition $S_k^{(0)}(n) \equiv S_k(n) = 1^k + 2^k + \cdots + n^k$, and $S_k^{(m)}(0) =0$. Recently, Karg{\i}n {\it et al.\/}, in a notable paper (arXiv:2008.00284 [math.NT]) obtained (among many other things) the following formula for $S_k^{(m)}(n)$ (see \cite[Equation (27)]{kargin})
\begin{equation}\label{kargin}
S_k^{(m)}(n) = \frac{1}{m!} \sum_{i=0}^{m} (-1)^i \genfrac{[}{]}{0pt}{}{m+n+1}{i+n+1}_{n+1} S_{k+i}(n),
\end{equation}
where $\genfrac{[}{]}{0pt}{}{m+r}{i+r}_{r}$ are the $r$-Stirling numbers of the first kind \cite{broder}. On the other hand, as shown in \cite{cere}, it turns out that (see \cite[Equation (2)]{cere})
\begin{equation}\label{cere}
S_k^{(m)}(n) = \frac{1}{m!} \sum_{i=0}^{m} (-1)^i q_{m,i}(n) S_{k+i}(n),
\end{equation}
where $q_{m,i}(n)$ is the polynomial in $n$ of degree $m-i$ given by (cf. \cite[Equation (16)]{cere})
\begin{equation*}
q_{m,i}(n) = \sum_{j=0}^{m-i} \binom{i+j}{i} \genfrac{[}{]}{0pt}{}{m+1}{i+j+1} n^j.
\end{equation*}

In view of \eqref{kargin} and \eqref{cere}, it follows that
\begin{equation}\label{intro}
\sum_{i=0}^{m} (-1)^i Q_{m,i}(n) S_{k+i}(n) =0,
\end{equation}
where $Q_{m,i}(n) = \genfrac{[}{]}{0pt}{}{m+n+1}{i+n+1}_{n+1} - q_{m,i}(n)$ is \emph{independent} of $k$. Moreover, it is to be noticed that \eqref{intro} holds irrespective of the value of $k$. Consequently, noting that the power sum polynomials $S_k(n), S_{k+1}(n),\ldots, S_{k+m}(n)$ are linearly independent, we are led to conclude that $Q_{m,i}(n) =0$ and then the following (nontrivial) relation must hold.

\begin{proposition}
For integers $m,i,n \geq 0$ and $m \geq i$, we have
\begin{equation}\label{result}
\genfrac{[}{]}{0pt}{}{m+n+1}{i+n+1}_{n+1} = \sum_{j=0}^{m-i} \binom{i+j}{i} \genfrac{[}{]}{0pt}{}{m+1}{i+j+1} n^j.
\end{equation}
\end{proposition}

In the next section we give an alternative derivation of this relation, showing that the formulas \eqref{kargin} and \eqref{cere} are in fact equivalent. In section 3, we show how the polynomial \eqref{result} relates to the so-called $r$-Stirling polynomials of the first kind. In section 4, we express the hyperharmonic polynomials and their successive derivatives in terms of the $r$-Stirling polynomials of the first kind. In section 5, we show the relationship between the (exponential) complete Bell polynomials and the $r$-Stirling numbers of the first kind. Finally, in section 6, we derive some identities involving the Bernoulli numbers and polynomials, the $r$-Stirling numbers of the first kind, the Stirling numbers of both kinds, and the harmonic numbers.

\section{Proof of the Proposition}

According to \cite[Equation (27)]{broder} (see also \cite[page 224]{mezo}), $\genfrac{[}{]}{0pt}{}{m+n+1}{i+n+1}_{n+1}$ can be expressed as
\begin{equation*}
\genfrac{[}{]}{0pt}{}{m+n+1}{i+n+1}_{n+1} = \sum_{j=i}^{m} \binom{m}{j} \genfrac{[}{]}{0pt}{}{j}{i} (n+1)^{\overline{m-j}},
\end{equation*}
where $n^{\overline{k}}$ denotes the rising factorial $n(n+1)\cdots (n+k-1)$. Changing the summation variable from $j$ to $t$, where $t = m-j$, results in
\begin{equation*}
\genfrac{[}{]}{0pt}{}{m+n+1}{i+n+1}_{n+1} = \sum_{t=0}^{m-i} \binom{m}{t} \genfrac{[}{]}{0pt}{}{m-t}{i} (n+1)^{\overline{t}}.
\end{equation*}
Now, we have that
\begin{equation*}
(n+1)^{\overline{t}} = \sum_{r=0}^{t} \genfrac{[}{]}{0pt}{}{t}{r} (n+1)^r = \sum_{r=0}^{t} \sum_{s=0}^{r}
\binom{r}{s} \genfrac{[}{]}{0pt}{}{t}{r} n^s,
\end{equation*}
and then
\begin{align*}
\genfrac{[}{]}{0pt}{}{m+n+1}{i+n+1}_{n+1} & = \sum_{t=0}^{m-i} \sum_{r=0}^{t} \sum_{s=0}^{r}
\binom{m}{t} \genfrac{[}{]}{0pt}{}{m-t}{i} \binom{r}{s} \genfrac{[}{]}{0pt}{}{t}{r} n^s \\
& = \sum_{t=0}^{m-i} \sum_{s=0}^{t} \sum_{r=s}^{t} \binom{r}{s} \genfrac{[}{]}{0pt}{}{t}{r}
\binom{m}{t} \genfrac{[}{]}{0pt}{}{m-t}{i} n^s .
\end{align*}
Using the well-known identity (see, e.g., \cite[Equation (6.16)]{graham}) $\sum_{r=s}^{t} \binom{r}{s} \genfrac{[}{]}{0pt}{}{t}{r} = \genfrac{[}{]}{0pt}{}{t+1}{s+1}$, the last equation becomes
\begin{align}\label{proof}
\genfrac{[}{]}{0pt}{}{m+n+1}{i+n+1}_{n+1} & = \sum_{t=0}^{m-i} \sum_{s=0}^{t}
\binom{m}{t} \genfrac{[}{]}{0pt}{}{m-t}{i} \genfrac{[}{]}{0pt}{}{t+1}{s+1} n^s \notag \\
& = \sum_{j=0}^{m-i} \left( \sum_{t=j}^{m-i} \binom{m}{t} \genfrac{[}{]}{0pt}{}{t+1}{j+1}
\genfrac{[}{]}{0pt}{}{m-t}{i} \right) n^j,
\end{align}
where we have renamed the variable $s$ as $j$.

Invoking the identity (see \cite[Equation (52)]{broder} and \cite[page 224]{mezo})
\begin{equation*}
\binom{j+i}{i} \genfrac{[}{]}{0pt}{}{m+r+s}{j+i+r+s}_{r+s} = \sum_{t=j}^{m-i} \binom{m}{t}
\genfrac{[}{]}{0pt}{}{t+r}{j+r}_{r} \genfrac{[}{]}{0pt}{}{m-t+s}{i+s}_{s},
\end{equation*}
and specializing to the case in which $r=1$ and $s=0$, it follows that
\begin{equation}\label{proof2}
\binom{i+j}{i} \genfrac{[}{]}{0pt}{}{m+1}{i+j+1} = \sum_{t=j}^{m-i} \binom{m}{t}
\genfrac{[}{]}{0pt}{}{t+1}{j+1} \genfrac{[}{]}{0pt}{}{m-t}{i},
\end{equation}
and then, combining \eqref{proof} and \eqref{proof2}, we get \eqref{result}.

\begin{remark}
Clearly, the constant term of the polynomial \eqref{result} is $\genfrac{[}{]}{0pt}{}{m+1}{i+1}$. Therefore, setting $n=0$ and $j=0$ in \eqref{proof}, we deduce the identity (cf. \cite[Equation (30)]{broder})
\begin{equation*}
\genfrac{[}{]}{0pt}{}{m+1}{i+1} = \sum_{t=0}^{m-i} t! \binom{m}{t} \genfrac{[}{]}{0pt}{}{m-t}{i}.
\end{equation*}
\end{remark}

\begin{remark}
When $k=0$, the hyper-sum $S_k^{(m)}(n)$ is given by $S_0^{(m)}(n) = \binom{n+m}{m+1}$. Hence, letting $k=0$ in \eqref{kargin} yields
\begin{equation}\label{recur}
\sum_{i=0}^{m} (-1)^i \genfrac{[}{]}{0pt}{}{m+n+1}{i+n+1}_{n+1} S_{i}(n) = m! \binom{n+m}{m+1}.
\end{equation}
In particular, for $n=1$ we find that
\begin{equation*}
\sum_{i=0}^{m} (-1)^i \genfrac{[}{]}{0pt}{}{m+2}{i+2}_{2} = m! .
\end{equation*}
\end{remark}

\section{$r$-Stirling polynomials of the first kind}

For integers $0 \leq i \leq m$, the $r$-Stirling polynomials of the first kind $R_{m,i}(x)$ are defined for arbitrary $x$ as (see \cite[Equation (56)]{broder} and \cite[Equation (5.3)]{carlitz})
\begin{equation*}
R_{m,i}(x) = \sum_{j=0}^{m-i} \binom{m}{j} \genfrac{[}{]}{0pt}{}{m-j}{i} x^{\overline{j}},
\end{equation*}
which reduces to $R_{m,i}(r) = \genfrac{[}{]}{0pt}{}{m+r}{i+r}_r$ when $r$ is a nonnegative integer. Following a procedure analogous to that used in the previous section, it can be shown that $R_{m,i}(x)$ can equally be written in the form (cf. \cite[Equation (5.2)]{carlitz})
\begin{equation}\label{carlitz}
R_{m,i}(x) = \sum_{j=0}^{m-i} \binom{i+j}{i} \genfrac{[}{]}{0pt}{}{m}{i+j} x^j.
\end{equation}

In view of \eqref{result}, it may be useful to define the following variant of \eqref{carlitz}
\begin{equation*}
\overline{R}_{m,i}(x) = \sum_{j=0}^{m-i} \binom{i+j}{i} \genfrac{[}{]}{0pt}{}{m+1}{i+j+1} x^j.
\end{equation*}
Observe that, for $x=-1$, we have
\begin{align*}
\overline{R}_{m,i}(-1) &= \sum_{j=0}^{m-i} (-1)^j \binom{i+j}{i} \genfrac{[}{]}{0pt}{}{m+1}{i+j+1} \\
& = \sum_{t=i}^{m} (-1)^{t-i} \binom{t}{i} \genfrac{[}{]}{0pt}{}{m+1}{t+1} = \genfrac{[}{]}{0pt}{}{m}{i},
\end{align*}
according to \cite[Equation (6.18)]{graham}. Moreover, $\overline{R}_{m,i}(0) = \genfrac{[}{]}{0pt}{}{m+1}{i+1}$, and, in general, $\overline{R}_{m,i}(r) = \genfrac{[}{]}{0pt}{}{m+r+1}{i+r+1}_{r+1} = R_{m,i}(r+1)$ for any integer $r \geq 0$. This in turn implies that, for arbitrary $x$,
\begin{equation}\label{myid}
\sum_{j=0}^{m-i} \binom{i+j}{i} \genfrac{[}{]}{0pt}{}{m+1}{i+j+1} x^j =
\sum_{j=0}^{m-i} \binom{i+j}{i} \genfrac{[}{]}{0pt}{}{m}{i+j} (x+1)^j.
\end{equation}
Incidentally, when $i=0$, \eqref{myid} reduces to
\begin{equation*}
\sum_{j=0}^{m} \genfrac{[}{]}{0pt}{}{m+1}{j+1} x^j = \sum_{j=0}^{m} \genfrac{[}{]}{0pt}{}{m}{j} (x+1)^j,
\end{equation*}
with each side of this equation being equal to the product $(x+1)(x+2)\cdots (x+m)$.

\section{Hyperharmonic polynomials and their derivatives}

The $n$th hyperharmonic number of order $r$, denoted by $H_n^{(r)}$, is defined by
\begin{equation*}
H_n^{(r)} = \left\{
              \begin{array}{ll}
                0 & \text{if}\,\, n \leq 0 \,\, \text{or}\,\, r <0,  \\
                \frac{1}{n}  &  \text{if}\,\, n >0 \,\, \text{and}\,\, r =0,\\
                \sum_{i=1}^n H_i^{(r-1)} &  \text{if}\,\, n,r \geq 1.
              \end{array}
            \right.
\end{equation*}
Note that $H_n^{(1)}$ is the ordinary harmonic number $H_n = 1 + \frac{1}{2} + \cdots + \frac{1}{n}$. As is well known (see \cite[Theorem 2]{benjamin}), $H_n^{(r)}$ can be given in terms of the $r$-Stirling number as follows: $H_n^{(r)} = \frac{1}{n!}\genfrac{[}{]}{0pt}{}{n+r}{r+1}_{r}$. Thus, employing \eqref{carlitz} allows us to express the $(j+1)$th hyperharmonic number of order $n$, $H_{j+1}^{(n)}$, in the form
\begin{equation*}
H_{j+1}^{(n)} = \frac{R_{j+1,1}(n)}{(j+1)!} = \frac{1}{(j+1)!} \sum_{i=0}^{j} (i+1) \genfrac{[}{]}{0pt}{}{j+1}{i+1} n^i,
\end{equation*}
giving $H_{j+1}^{(n)}$ as a polynomial in $n$ of degree $j$. Let us observe that, when $n=1$, we recover the well-known identity $H_j = \frac{1}{j!} \sum_{i=1}^{j} i \genfrac{[}{]}{0pt}{}{j}{i}$. Of course, $H_{j+1}^{(n)}$ can be naturally extended to a polynomial $H_{j+1}^{(x)}$ in which $x$ takes any arbitrary value as follows
\begin{equation}\label{hyper}
H_{j+1}^{(x)} = \frac{R_{j+1,1}(x)}{(j+1)!} = \frac{1}{(j+1)!} \sum_{i=0}^{j} (i+1) \genfrac{[}{]}{0pt}{}{j+1}{i+1} x^i.
\end{equation}
Likewise, we have
\begin{equation}\label{hyper2}
H_{j+1}^{(x+1)} = \frac{\overline{R}_{j+1,1}(x)}{(j+1)!} = \frac{1}{(j+1)!} \sum_{i=0}^{j} (i+1)
\genfrac{[}{]}{0pt}{}{j+2}{i+2} x^i,
\end{equation}
giving $H_{j+1}^{(x+1)}$ as a polynomial in $x$ of degree $j$.

On the other hand, from \cite[Theorem 28]{broder}, we know that $R_{m,i}(x) = \frac{1}{i!} \frac{d^{i}}{d x^{i}} x^{\overline{m}}$. Therefore, it follows that
\begin{equation*}
H_{j+1}^{(x)} = \frac{1}{(j+1)!} \frac{d}{dx} x^{\overline{j+1}} = \frac{d}{dx}
\binom{x+j}{j+1},
\end{equation*}
in accordance with \cite[Corollary 1]{cere2}. Moreover, the $i$th derivative of $H_{j+1}^{(x)}$ with respect to $x$ is given by
\begin{equation*}
\frac{d^{i}}{d x^i}H_{j+1}^{(x)} = \frac{1}{(j+1)!} \frac{d^{i+1}}{d x^{i+1}}x^{\overline{j+1}}
= \frac{(i+1)!}{(j+1)!} R_{j+1,i+1}(x),
\end{equation*}
and then
\begin{equation}\label{der1}
\frac{d^{i}}{d x^i}H_{j+1}^{(x)} = \frac{(i+1)!}{(j+1)!} \sum_{t=0}^{j-i} \binom{i+t+1}{i+1}
\genfrac{[}{]}{0pt}{}{j+1}{i+t+1} x^t.
\end{equation}
Likewise, we have
\begin{equation*}
\frac{d^{i}}{d x^i}H_{j+1}^{(x+1)} =  \frac{(i+1)!}{(j+1)!} R_{j+1,i+1}(x+1)
= \frac{(i+1)!}{(j+1)!} \overline{R}_{j+1,i+1}(x),
\end{equation*}
and then
\begin{equation}\label{der2}
\frac{d^{i}}{d x^i}H_{j+1}^{(x+1)} = \frac{(i+1)!}{(j+1)!} \sum_{t=0}^{j-i} \binom{i+t+1}{i+1}
\genfrac{[}{]}{0pt}{}{j+2}{i+t+2} x^t.
\end{equation}
Note that, as it must be, \eqref{der1} and \eqref{der2} reduce respectively to \eqref{hyper} and \eqref{hyper2} when $i=0$.
\pagebreak

Furthermore, for nonnegative integer $r$, from \eqref{der1} and \eqref{der2} it follows that
\begin{equation*}
\left. \frac{d^{i}}{d x^i}H_{j+1}^{(x)}\right|_{x =r} = \frac{(i+1)!}{(j+1)!} \genfrac{[}{]}{0pt}{}{j+r+1}{i+r+1}_r,
\end{equation*}
and
\begin{equation*}
\left. \frac{d^{i}}{d x^i}H_{j+1}^{(x+1)}\right|_{x =r} = \frac{(i+1)!}{(j+1)!} \genfrac{[}{]}{0pt}{}{j+r+2}{i+r+2}_{r+1},
\end{equation*}
In particular, for $r=0$ we obtain
\begin{equation*}
\left. \frac{d^{i}}{d x^i}H_{j+1}^{(x)}\right|_{x =0} = \frac{(i+1)!}{(j+1)!} \genfrac{[}{]}{0pt}{}{j+1}{i+1},
\end{equation*}
and
\begin{equation*}
\left. \frac{d^{i}}{d x^i}H_{j+1}^{(x+1)}\right|_{x =0} = \frac{(i+1)!}{(j+1)!} \genfrac{[}{]}{0pt}{}{j+2}{i+2}.
\end{equation*}

\begin{remark}
According to \cite[Theorem 2]{cere2}, the hyperharmonic polynomial $H_{j+1}^{(x+1)}$ can be expressed as
\begin{equation*}
H_{j+1}^{(x+1)} = \sum_{t=0}^{j} \frac{1}{j+1-t} \binom{x+t}{t}.
\end{equation*}
Using this representation in combination with \eqref{der2} gives
\begin{equation*}
\sum_{t=0}^{j} \frac{1}{j+1-t} \, \frac{d^{i}}{d x^i} \binom{x+t}{t} = \frac{(i+1)!}{(j+1)!} \sum_{t=0}^{j-i}
\binom{i+t+1}{i+1} \genfrac{[}{]}{0pt}{}{j+2}{i+t+2} x^t.
\end{equation*}
As a consequence, we have
\begin{equation*}
\sum_{t=0}^{j} \frac{1}{j+1-t} \left. \frac{d^{i}}{d x^i} \binom{x+t}{t}\right|_{x=r}
= \frac{(i+1)!}{(j+1)!} \genfrac{[}{]}{0pt}{}{j+r+2}{i+r+2}_{r+1},
\end{equation*}
or, equivalently,
\begin{equation}\label{wuyu}
\sum_{t=0}^{j} \frac{1}{t! \, (j+1-t)} \genfrac{[}{]}{0pt}{}{t+r}{i+r}_{r}
= \frac{i+1}{(j+1)!} \genfrac{[}{]}{0pt}{}{j+r+1}{i+r+1}_{r},
\end{equation}
which holds for any integers $0 \leq i \leq  j$ and $r \geq 0$. In particular, substituting $i=r=1$ in \eqref{wuyu} yields the identity
\begin{equation*}
\sum_{t=0}^{j} \frac{H_t}{j+1-t} = \frac{2}{(j+1)!} \genfrac{[}{]}{0pt}{}{j+2}{3} = H_{j+1}^2 - H_{j+1}^{[2]},
\end{equation*}
which can also be found in \cite[p. 544]{kargin2}, and where the notation $H_{j}^{[2]}$ means $\sum_{t=1}^{j} 1/t^2$.
\end{remark}

\section{$r$-Stirling numbers of the first kind and complete Bell polynomials}

From the previous section we know that, for $j \geq 1$,
\begin{equation}\label{myres}
\left. \frac{d^{i}}{d x^i}H_{j}^{(x+1)}\right|_{x =r} = \frac{(i+1)!}{j!} \genfrac{[}{]}{0pt}{}{j+r+1}{i+r+2}_{r+1}.
\end{equation}
In what follows we are going to show that \eqref{myres} is consistent with the result found in \cite[Equation (19)]{kargin}. To this end, we first introduce the concepts behind the said equation \cite[Equation (19)]{kargin}.

Let $Y_n(x_1,x_2,\ldots,x_n)$ be the $n$th (exponential) complete Bell polynomial defined by $Y_0 =1$ and (\cite[page 134]{comtet})
\begin{equation*}
\exp \left( \sum_{j=1}^{\infty} x_j \frac{t^j}{j!} \right) = 1 + \sum_{p=1}^{\infty} Y_p(x_1,x_2,\ldots,x_p)
\frac{t^p}{p!}.
\end{equation*}
Now, for integers $r \geq 0$ and $j,k \geq 1$, let us define $H(j,k;r)$ as
\begin{equation*}
H(j,k;r) = \sum_{t=1}^{j} \frac{1}{(t+r)^k}.
\end{equation*}
Note that $H(j,1;0) = H_j$. Moreover, following Spie{\ss} (see \cite[Equation (6)]{spieb}), we introduce the numbers $P(i,j+r,r)$ as
\begin{equation*}
P(i,j+r,r) = P_i \big( H(j,1;r), H(j,2;r), \ldots, H(j,i;r) \big),
\end{equation*}
where the polynomial $P_i(x_1,x_2,\ldots, x_i)$ is defined by
\begin{equation*}
P_i(x_1,x_2,\ldots, x_i) = (-1)^i Y_i (-0! x_1, -1!x_2, \ldots, -(i-1)! x_i),
\end{equation*}
or, equivalently,
\begin{equation*}
P_i(x_1,x_2,\ldots, x_i) = Y_i (0! x_1, -1!x_2, \ldots, (-1)^{i-1}(i-1)! x_i).
\end{equation*}
The first five polynomials $P_i(x_1,x_2,\ldots, x_i)$ are given by
\begin{align*}
  & P_1(x_1) = x_1, \\
  & P_2(x_1,x_2) = x_1^2 - x_2, \\
  & P_3(x_1,x_2,x_3) = x_1^3 -3x_1 x_2 + 2x_3, \\
  & P_4(x_1,x_2,x_3,x_4) = x_1^4 -6x_1^2 x_2 + 8x_1 x_3 + 3x_2^2 -6x_4, \\
  & P_4(x_1,x_2,x_3,x_4,x_5) = x_1^5 -10x_1^3 x_2 + 20x_1^2 x_3 + 15x_1 x_2^2 -30x_1 x_4 -20x_2 x_3 +24x_5.
\end{align*}

With these ingredients at hand, we are ready to state the result established in \cite[Equation (19)]{kargin}, namely
\begin{equation}\label{hisres}
\left. \frac{d^{i}}{d x^i}H_{j}^{(x+1)}\right|_{x =r} = \binom{j+r}{r} P(i+1,j+r,r).
\end{equation}
Next we show that the right-hand sides of \eqref{myres} and \eqref{hisres} are the same. This is a direct consequence of the following theorem set forth by K\"{o}lbig in \cite[Theorem]{kolbig}.

\begin{theorem}[K\"{o}lbig, 1994]
Let $\alpha \in \mathbb{R}$ with $\alpha \neq -1, -2, \ldots, -j$, and
\begin{equation*}
H(j,k;\alpha) = \sum_{t=1}^{j} \frac{1}{(t+\alpha)^{k}},
\end{equation*}
for integers $j,k \geq 1$. Then we have
\begin{equation*}
P_q \big( H(j,1;\alpha), H(j,2;\alpha), \ldots, H(j,q;\alpha) \big) =
\left\{
              \begin{array}{ll}
                \displaystyle{\frac{q!}{(1+\alpha)^{\overline{j}}}} S(j,q;\alpha), & q \leq j,  \\
                0, &  q >j ,
              \end{array}
            \right.
\end{equation*}
where
\begin{equation*}
S(j,q;\alpha) = \sum_{t=q}^{j} \binom{t}{q} \genfrac{[}{]}{0pt}{}{j}{t} (1+ \alpha)^{t-q}.
\end{equation*}
\end{theorem}

Hence, taking $\alpha =r$ in K\"{o}lbig's theorem (with integer $r \geq 0$), it follows that
\begin{align*}
\binom{j+r}{r} P(i+1,j+r,r) & = \binom{j+r}{r} P_{i+1} \big( H(j,1;r), H(j,2;r), \ldots, H(j,i+1;r) \big) \\
& = \frac{(i+1)!}{(1+r)^{\overline{j}}} \binom{j+r}{r} S(j,i+1;r),
\end{align*}
where we assume that $j \geq 1$ and $j \geq i+1$, with $i =0,1,2,\ldots\,$. Furthermore, it turns out that $\binom{j+r}{r} = \frac{1}{j!} (1+r)^{\overline{j}}$, and then
\begin{equation*}
\binom{j+r}{r} P(i+1,j+r,r) = \frac{(i+1)!}{j!} S(j,i+1;r),
\end{equation*}
where
\begin{equation*}
S(j,i+1;r) = \sum_{t=i+1}^{j} \binom{t}{i+1} \genfrac{[}{]}{0pt}{}{j}{t} (1+r)^{t-i-1}.
\end{equation*}
On the other hand, we have
\begin{align*}
\genfrac{[}{]}{0pt}{}{j+r+1}{i+r+2}_{r+1} & = R_{j,i+1}(r+1) = \sum_{s=0}^{j-i-1} \binom{i+s+1}{i+1}
\genfrac{[}{]}{0pt}{}{j}{i+s+1} (r+1)^{s} \\
& = \sum_{t=i+1}^{j} \binom{t}{i+1} \genfrac{[}{]}{0pt}{}{j}{t} (r+1)^{t-i-1}.
\end{align*}
Therefore, we conclude that $S(j,i+1;r) = \genfrac{[}{]}{0pt}{}{j+r+1}{i+r+2}_{r+1}$, and we are done.

\begin{remark}
It should be remarked that the identity \eqref{hisres} was already established in \cite[Equation (4.2)]{wang} in the equivalent form
\begin{equation*}
\left. \frac{d^{i}}{d x^i} x^{\overline{j}} \right|_{x =r+1} = j! \binom{j+r}{r} P(i,j+r,r).
\end{equation*}
Furthermore, in \cite[p. 270]{wang2} we find that
\begin{equation*}
\left. \frac{d^{i}}{d x^i} \binom{x+n}{m} \right|_{x =0} = \binom{n}{m} P(i,n,n-m), \quad n \geq m >0,
\end{equation*}
which, as can be easily verified, is also equivalent to \eqref{hisres}.
\end{remark}

\begin{remark}
For any integer $r \geq 1$, when $i=0$ and $r \to r-1$, \eqref{hisres} reduces to
\begin{align*}
\left. H_{j}^{(x+1)}\right|_{x =r-1} & = H_j^{(r)} = \binom{j+r-1}{r-1} H(j,1;r-1) \\
& = \binom{j+r-1}{r-1} \sum_{t=1}^{j} \frac{1}{t+r-1} \\
& = \binom{j+r-1}{r-1} \big( H_{j+r-1} - H_{r-1} \big),
\end{align*}
thus retrieving the well-known formula connecting hyperharmonic numbers with ordinary harmonic numbers put forward by Conway and Guy \cite[p. 258]{conway}.
\end{remark}

Equating the right-hand sides of \eqref{myres} and \eqref{hisres}, and letting $i \to i-1$, we obtain the identity
\begin{equation}\label{myres2}
\binom{j+r}{r} P(i,j+r,r) = \frac{i!}{j!} \genfrac{[}{]}{0pt}{}{j+r+1}{i+r+1}_{r+1},
\end{equation}
connecting the complete Bell polynomials and the $r$-Stirling numbers of the first kind (cf. \cite[Equations (4.2) and (4.36)]{wang}). When $r=0$, we recover the well-known result (see \cite[Equation (7b), p. 217]{comtet})
\begin{equation*}
\genfrac{[}{]}{0pt}{}{j+1}{i+1} = \frac{j!}{i!} P(i,j,0) = \frac{j!}{i!} P_i \big( H_j^{[1]}, H_j^{[2]}, \ldots, H_j^{[i]} \big),
\end{equation*}
where $H_j^{[i]} \equiv H(j,i;0) = \sum_{t=1}^{j} 1/t^i$, and $H_j^{[1]} = H_j$.

\enlargethispage{4mm}
\begin{remark}
The generating function of the numbers $\binom{j+r}{r} P(i,j+r,r)$ is given by (see, e.g., \cite[Equation (1.6)]{wang})
\begin{equation*}
\sum_{j=i}^{\infty} \binom{j+r}{r} P(i,j+r,r) t^j = \frac{( -\ln(1-t))^i}{(1-t)^{r+1}}.
\end{equation*}
Therefore, using \eqref{myres2} and letting $r \to r-1$, we readily obtain the exponential generating function of the $r$-Stirling numbers of the first kind, namely
\begin{equation*}
\sum_{j=i}^{\infty} \genfrac{[}{]}{0pt}{}{j+r}{i+r}_{r} \frac{t^j}{j!} = \frac{1}{i!}
\frac{( -\ln(1-t))^i}{(1-t)^{r}}.
\end{equation*}
In particular, for $r=0$ and $i \geq 1$, we have
\begin{equation*}
( -\ln(1-t))^i = \sum_{j=i}^{\infty} \frac{i!}{j!} \genfrac{[}{]}{0pt}{}{j}{i} t^j = \sum_{j=i}^{\infty} \frac{i}{j} P(i-1,j-1,0) t^j,
\end{equation*}
in agreement with \cite[Theorem 9]{spieb}. Note that setting $i=1$ in last equation yields the Maclaurin series of the natural logarithm
\begin{equation*}
\ln(1-t) = - \sum_{j=1}^{\infty} \frac{t^j}{j} = -t - \frac{t^2}{2} - \frac{t^3}{3} - \ldots \, .
\end{equation*}
\end{remark}

\begin{remark}
According to \cite[Theorem 16]{spieb}, it turns out that, for $m,r \geq 0$,
\begin{equation*}
\sum_{k=0}^{m} \frac{P(r,k,0)}{k+1} = \frac{P(r+1,m+1,0)}{r+1}.
\end{equation*}
Recalling that $P(r,k,0) = \frac{r!}{k!}  \genfrac{[}{]}{0pt}{}{k+1}{r+1}$, it is easily seen that the above identity is equivalent to (cf. \cite[Equation (6.21)]{graham})
\begin{equation*}
\genfrac{[}{]}{0pt}{}{m+1}{r+1} = m! \sum_{k=0}^{m} \frac{1}{k!} \genfrac{[}{]}{0pt}{}{k}{r}.
\end{equation*}
\end{remark}

\begin{remark}
By virtue of \eqref{myres2}, we can reformulate \cite[Theorem 6]{kargin} as follows:
\begin{equation}\label{fin1}
\sum_{k=l}^{n} \genfrac{[}{]}{0pt}{}{n+r}{k+r}_{r} \binom{k}{l} B_{k-l}(q) =
\frac{l+1}{n+1} \genfrac{[}{]}{0pt}{}{n+q+r}{l+q+r}_{q+r-1},
\end{equation}
which may also be written as
\begin{equation}\label{fin2}
\sum_{k=l}^{n} \genfrac{[}{]}{0pt}{}{n+r+1}{k+r+1}_{r+1} \binom{k}{l} B_{k-l}(q) =
\frac{l+1}{n+1} \genfrac{[}{]}{0pt}{}{n+q+r+1}{l+q+r+1}_{q+r},
\end{equation}
or else,
\begin{equation}\label{fin3}
\sum_{k=l}^{n} \genfrac{[}{]}{0pt}{}{n+r}{k+r}_r \binom{k}{l} B_{k-l}(q+1) =
\frac{l+1}{n+1} \genfrac{[}{]}{0pt}{}{n+q+r+1}{l+q+r+1}_{q+r},
\end{equation}
for nonnegative integers $l,q,r$, and $n \geq l$. As an example, putting $l=2$, $q=0$, and $r=1$ in \eqref{fin3}, we get
\begin{align*}
\sum_{k=2}^{n} (-1)^k \genfrac{[}{]}{0pt}{}{n+1}{k+1} k(k-1) B_{k-2} & =
\frac{6}{n+1} \genfrac{[}{]}{0pt}{}{n+2}{4} \\
& = n! \big( H_{n+1}^3 -3 H_{n+1} H_{n+1}^{[2]} + 2H_{n+1}^{[3]} \big),
\end{align*}
in accordance with the particular identity found in \cite[page 8]{kargin}. By the way, from \eqref{fin2} and \eqref{fin3}, we find that
\begin{equation*}
\sum_{k=l}^{n} \genfrac{[}{]}{0pt}{}{n+r+1}{k+r+1}_{r+1} \binom{k}{l} B_{k-l}(x) =
\sum_{k=l}^{n} \genfrac{[}{]}{0pt}{}{n+r}{k+r}_r \binom{k}{l} B_{k-l}(x+1),
\end{equation*}
which holds for arbitrary $x$. In particular, for $r=0$, and renaming the indices $k-l \to j$, $l \to i$, and $n \to m$, the last identity becomes
\begin{equation*}
\sum_{j=0}^{m-i} \binom{i+j}{i} \genfrac{[}{]}{0pt}{}{m+1}{i+j+1} B_{j}(x) =
\sum_{j=0}^{m-i} \binom{i+j}{i} \genfrac{[}{]}{0pt}{}{m}{i+j} B_{j}(x+1),
\end{equation*}
which may be compared with \eqref{myid}.

Moreover, \eqref{fin1} can be generalized as follows:
\begin{equation}\label{fin4}
\sum_{k=l}^{n} \genfrac{[}{]}{0pt}{}{n+r}{k+r}_{r} \binom{k}{l} B_{k-l}(x) =
\frac{l+1}{n+1} \sum_{k=l}^{n} \binom{k+1}{l+1} \genfrac{[}{]}{0pt}{}{n+1}{k+1} (x+r-1)^{k-l},
\end{equation}
which holds for arbitrary $x$. For $l=0$, \eqref{fin4} can be compactly written as
\begin{equation}\label{fin5}
\sum_{k=0}^{n} \genfrac{[}{]}{0pt}{}{n+r}{k+r}_{r} B_{k}(x) = n! H_{n+1}^{(x+r-1)}.
\end{equation}
Conversely, \eqref{fin4} can be obtained by performing the $l$th derivative with respect to $x$ of both sides of \eqref{fin5}. For $x = -r$, \eqref{fin5} reads
\begin{equation*}
\sum_{k=0}^{n} \genfrac{[}{]}{0pt}{}{n+r}{k+r}_{r} B_{k}(-r) = n! H_{n+1}^{(-1)}.
\end{equation*}
Now, according to \cite[Equation (42)]{cere2}, we have
\begin{equation*}
H_{n+1}^{(-1)} =
\left\{
  \begin{array}{ll}
   -\dfrac{1}{n(n+1)}, &  n \geq 1, \\
    1, & n=0,
  \end{array}
\right.
\end{equation*}
and then we get the identity
\begin{equation*}
\sum_{k=0}^{n} \genfrac{[}{]}{0pt}{}{n+r}{k+r}_{r} B_{k}(-r) = -\frac{n!}{n(n+1)},
\end{equation*}
which holds for any integers $r \geq 0$ and $n \geq 1$. In particular, for $r=0$ we obtain
\begin{equation*}
\sum_{k=0}^{n} \genfrac{[}{]}{0pt}{}{n}{k} B_{k} = -\frac{(n-1)!}{n+1}, \quad n \geq 1.
\end{equation*}
\end{remark}

\begin{remark}
Using \eqref{myres2} into \eqref{recur} leads to the relation
\begin{equation*}
\sum_{i=0}^m \frac{(-1)^i}{i!} P(i,m+n,n) S_i(n) = \frac{n}{m+1},
\end{equation*}
and, in particular,
\begin{equation*}
\sum_{i=0}^m \frac{(-1)^i}{i!} P(i,m+1,1) = \frac{1}{m+1}.
\end{equation*}
\end{remark}

\section{Further connections with $r$-Stirling numbers of the first kind}

In this section we will provide further identities involving the $r$-Stirling numbers of the first kind, the Bernoulli numbers and polynomials, the ordinary Stirling numbers of the first and second kinds, and the harmonic numbers.

Firstly, from \cite[Equation (69)]{cere2}, it follows that, for $k \geq 0$, the Bernoulli polynomials $B_k(x)$ can be expressed in terms of $H_{j+1}^{(x)}$ as
\begin{equation}\label{ber}
B_k(x) = (-1)^k \sum_{j=0}^{k} (-1)^{j} j! \genfrac{\{}{\}}{0pt}{}{k+1}{j+1} H_{j+1}^{(x)},
\end{equation}
where $\genfrac{\{}{\}}{0pt}{}{k}{j}$ are the Stirling numbers of the second kind. Hence, combining \eqref{hyper} and \eqref{ber} gives
\begin{equation*}
B_k(x) = (-1)^k \sum_{i=0}^{k} \sum_{j=i}^{k} (-1)^j \, \frac{i+1}{j+1} \genfrac{\{}{\}}{0pt}{}{k+1}{j+1}
\genfrac{[}{]}{0pt}{}{j+1}{i+1} x^i,
\end{equation*}
and, for $x=0$,
\begin{equation*}
B_k = (-1)^k \sum_{j=0}^{k} (-1)^j \frac{j!}{j+1} \genfrac{\{}{\}}{0pt}{}{k+1}{j+1}.
\end{equation*}
Furthermore, taking successive derivatives of \eqref{ber} we find that
\begin{align*}
\frac{d^{i}}{d x^i}B_k(x) & = i! \binom{k}{i}B_{k-i}(x) \\
& = (-1)^k (i+1)! \sum_{j=i}^{k} \frac{(-1)^j}{j+1} \genfrac{\{}{\}}{0pt}{}{k+1}{j+1} R_{j+1,i+1}(x),
\end{align*}
where we have used that $R_{j,i}(x) =0$ for $i >j$. When $r$ is a nonnegative integer, we obtain
\begin{equation*}
\left. \frac{d^{i}}{d x^i}B_k(x) \right|_{x=r} = (-1)^k (i+1)! \sum_{j=i}^{k} \frac{(-1)^j}{j+1} \genfrac{\{}{\}}{0pt}{}{k+1}{j+1}
\genfrac{[}{]}{0pt}{}{j+r+1}{i+r+1}_r ,
\end{equation*}
and, in particular,
\begin{equation*}
B_k(r) = (-1)^k \sum_{j=0}^{k} \frac{(-1)^j}{j+1} \genfrac{\{}{\}}{0pt}{}{k+1}{j+1} \genfrac{[}{]}{0pt}{}{j+r+1}{r+1}_r.
\end{equation*}

On the other hand,  Spie{\ss} \cite{spieb} and Wang \cite{wang} derived a number of identities involving the numbers $P(i,j+r,r)$. Next, making use of \eqref{myres2}, we rewrite some of these identities in terms of the $r$-Stirling numbers $\genfrac{[}{]}{0pt}{}{j+r+1}{i+r+1}_{r+1}$.

\begin{itemize}

\item
From \cite[Theorem 10]{spieb} and \eqref{myres2}, we get
\begin{equation*}
\sum_{k=0}^{m} \frac{(-1)^k}{k!} \binom{r+1}{m-k} \genfrac{[}{]}{0pt}{}{k+r+1}{i+r+1}_{r+1}
= \frac{(-1)^m}{m!} \genfrac{[}{]}{0pt}{}{m}{i}.
\end{equation*}
When $r=0$, this identity yields the well-known recurrence relation $\genfrac{[}{]}{0pt}{}{m+1}{i+1} = \genfrac{[}{]}{0pt}{}{m}{i} + m \genfrac{[}{]}{0pt}{}{m}{i+1}$.

\item
From \cite[Theorem 13]{spieb} and \eqref{myres2}, we get
\begin{equation*}
\sum_{k=r}^{m+1-s} \frac{r! s!}{k! \, (m+1-k)!} \genfrac{[}{]}{0pt}{}{k}{r} \genfrac{[}{]}{0pt}{}{m+1-k}{s}
= \frac{(r+s)!}{(m+1)!} \genfrac{[}{]}{0pt}{}{m+1}{r+s}.
\end{equation*}
In particular, for $s=1$ we obtain
\begin{equation*}
\sum_{k=r}^{m} \frac{1}{k! \, (m+1-k)} \genfrac{[}{]}{0pt}{}{k}{r} = \frac{(r+1)}{(m+1)!} \genfrac{[}{]}{0pt}{}{m+1}{r+1},
\end{equation*}
which is just the identity in \eqref{wuyu} evaluated for $r=0$.

\item
From \cite[Theorem 15]{spieb} and \eqref{myres2}, and after a few simple manipulation, we get
\begin{equation}\label{spi}
q^{r} \sum_{k=1}^{m} \frac{1}{(k+q)!} \genfrac{[}{]}{0pt}{}{k}{r} = \frac{1}{q!} - \frac{1}{(m+q)!}
\sum_{j=1}^{r} q^{j-1} \genfrac{[}{]}{0pt}{}{m+1}{j},
\end{equation}
which is valid for any integers $q \geq 1$ and $1 \leq r \leq m$. Setting $r=1,2$, and 3 in \eqref{spi} gives
\begin{align*}
\sum_{k=1}^m \frac{1}{k(k+1) \cdots (k+q)} = & \frac{1}{q \cdot q!} - \frac{1}{q (m+1) \cdots (m+q)}, \\
\sum_{k=1}^m \frac{H_{k-1}}{k(k+1) \cdots (k+q)} = & \frac{1}{q^2 \cdot q!} - \frac{1 + qH_m}{q^2 (m+1) \cdots (m+q)}, \\
\intertext{and}
\sum_{k=1}^m \frac{H_{k-1}^2 - H_{k-1}^{[2]}}{k(k+1) \cdots (k+q)} = & \frac{2}{q^3 \cdot q!} - \frac{2 + 2q H_m +
q^2 \big(H_m^{2} - H_m^{[2]}\big)}{q^3 (m+1) \cdots (m+q)},
\end{align*}
respectively. The above three identities are to be compared with the corresponding Examples 1, 2, and 3 previously obtained in \cite[p. 849]{spieb}. Moreover, it is worth pointing out that, for the case in which $r =m$, \eqref{spi} yields the horizontal generating function for the Stirling numbers $\genfrac{[}{]}{0pt}{}{m+1}{j+1}$, namely
\begin{equation*}
\sum_{j=0}^{m} \genfrac{[}{]}{0pt}{}{m+1}{j+1} q^j = (q+1)(q+2) \cdots (q+m),
\end{equation*}
which holds for arbitrary $q$.

\item
From \cite[Equation (3.4)]{wang} and \eqref{myres2}, we get
\begin{equation*}
\sum_{k=0}^{m} \frac{(-1)^k}{k!} \binom{m}{k} \binom{k+r}{r}^{-1} \genfrac{[}{]}{0pt}{}{k+r+1}{i+r+1}_{r+1}
= \frac{(-1)^i}{m!} \binom{m+r}{r}^{-1} \genfrac{[}{]}{0pt}{}{m}{i}.
\end{equation*}
Letting $r=i=1$ in the above relation leads to
\begin{equation*}
\sum_{k=0}^{m} (-1)^k \binom{m}{k} H_{k+1} = -\frac{1}{m(m+1)},
\end{equation*}
which may be compared with the more commonly known identity (see, e.g., \cite[Equation (3.2)]{wang2}) $\sum_{k=0}^{m} (-1)^k \binom{m}{k} H_{k} = -\frac{1}{m}$, where $m \geq 1$.

\item
From \cite[Equation (3.22)]{wang} and \eqref{myres2}, we get
\begin{equation*}
\genfrac{[}{]}{0pt}{}{m+r+1}{i+r+1}_{r+1} = m! \sum_{k=0}^{m} \frac{1}{k!} \binom{m-k+r}{r}
\genfrac{[}{]}{0pt}{}{k}{i} = m! \sum_{k=0}^{m} (-1)^{k-i} \frac{1}{k!} \binom{m+r}{k+r} \genfrac{[}{]}{0pt}{}{k}{i}.
\end{equation*}
In particular, for $r=0$ we have (cf. \cite[Equation (3.25)]{wang})
\begin{equation*}
\genfrac{[}{]}{0pt}{}{m+1}{i+1} = m! \sum_{k=0}^{m} \frac{1}{k!} \genfrac{[}{]}{0pt}{}{k}{i}
= m! \sum_{k=0}^{m} (-1)^{k-i} \frac{1}{k!} \binom{m}{k} \genfrac{[}{]}{0pt}{}{k}{i}.
\end{equation*}

\item
From \cite[Equations (3.32) and (3.33)]{wang} and \eqref{myres2}, we get
\begin{equation}\label{wa1}
\sum_{k=0}^{m} \frac{1}{k!} \binom{r+m-k-1}{m-k} \genfrac{[}{]}{0pt}{}{k+s+1}{i+s+1}_{s+1}
= \frac{1}{m!} \genfrac{[}{]}{0pt}{}{m+r+s+1}{i+r+s+1}_{r+s+1},
\end{equation}
and
\begin{equation}\label{wa2}
\genfrac{[}{]}{0pt}{}{m+r+1}{i+r+1}_{r+1} = m! \sum_{k=0}^{m} \frac{1}{k!} \binom{r+m-k-1}{m-k}
\genfrac{[}{]}{0pt}{}{k+1}{i+1},
\end{equation}
respectively. Notice that, from \eqref{wa1} and \eqref{wa2}, one quickly obtains
\begin{equation*}
\sum_{k=0}^{m} \frac{1}{k!} \binom{r+m-k-1}{m-k} \genfrac{[}{]}{0pt}{}{k+s+1}{i+s+1}_{s+1}
= \sum_{k=0}^{m} \frac{1}{k!} \binom{r+s+m-k-1}{m-k} \genfrac{[}{]}{0pt}{}{k+1}{i+1}.
\end{equation*}
Furthermore, regarding \eqref{wa2}, for $r=1$ it reads
\begin{equation*}
\genfrac{[}{]}{0pt}{}{m+2}{i+2}_{2} = m! \sum_{k=0}^{m} \frac{1}{k!} \genfrac{[}{]}{0pt}{}{k+1}{i+1},
\end{equation*}
or, equivalently,
\begin{equation*}
P_i \big( H(m,1;1), H(m,2;1), \ldots, H(m,i;1) \big) = \frac{i!}{m+1} \sum_{k=0}^{m} \frac{1}{k!} \genfrac{[}{]}{0pt}{}{k+1}{i+1}.
\end{equation*}
In particular, when $i=1$ we recover the well-known identity
\begin{equation*}
\sum_{k=0}^{m} H_k = (m+1)H_m -m.
\end{equation*}

\item
From \cite[Equations (4.1) and (4.3)]{wang} and \eqref{myres2}, we get
\begin{equation}\label{it1}
\sum_{k=0}^{m} (-1)^{m-k} \genfrac{\{}{\}}{0pt}{}{m}{k} \genfrac{[}{]}{0pt}{}{k+r+1}{i+r+1}_{r+1}
= \binom{m}{i} (r+1)^{m-i},
\end{equation}
and
\begin{equation}\label{it2}
\sum_{k=i}^{m} (-r-1)^{k-i} \binom{k}{i} \genfrac{[}{]}{0pt}{}{m+r+1}{k+r+1}_{r+1}
= \genfrac{[}{]}{0pt}{}{m}{i},
\end{equation}
respectively. When $r=0$, \eqref{it1} and \eqref{it2} reduce to
\begin{equation*}
\sum_{k=0}^{m} (-1)^{m-k} \genfrac{\{}{\}}{0pt}{}{m}{k} \genfrac{[}{]}{0pt}{}{k+1}{i+1}
= \binom{m}{i},
\end{equation*}
and
\begin{equation*}
\sum_{k=i}^{m} (-1)^{k-i} \binom{k}{i} \genfrac{[}{]}{0pt}{}{m+1}{k+1}
= \genfrac{[}{]}{0pt}{}{m}{i},
\end{equation*}
respectively. Furthermore, when $r=i=1$, from \eqref{it1} we obtain
\begin{equation*}
\sum_{k=0}^{m} (-1)^{m-k} (k+1)! \genfrac{\{}{\}}{0pt}{}{m}{k} \big( H_{k+1} -1 \big) = m 2^{m-1}.
\end{equation*}
On the other hand, putting $i=1$ in \cite[Equation (4.5)]{wang} yields
\begin{equation*}
\sum_{k=0}^{m} (-1)^{m-k} (k+1)! \genfrac{\{}{\}}{0pt}{}{m}{k} = 2^{m}.
\end{equation*}
Therefore, combining the last two identities, we get
\begin{equation*}
\sum_{k=0}^{m} (-1)^{m-k} (k+1)! \genfrac{\{}{\}}{0pt}{}{m}{k} H_{k+1} = (m+2) \, 2^{m-1},
\end{equation*}
which may be compared with \cite[Equation (4.10)]{wang}.

\item
From \cite[Equation (4.40)]{wang} and \eqref{myres2}, we get
\begin{equation*}
\sum_{k=0}^{m} (-1)^{k-i} i! \genfrac{\{}{\}}{0pt}{}{k}{i} \genfrac{[}{]}{0pt}{}{m+r+1}{k+r+1}_{r+1}
= m! \binom{r+m-i}{m-i},
\end{equation*}
which, for $r=0$, reduces to
\begin{equation*}
\sum_{k=0}^{m} (-1)^{k-i} i! \genfrac{\{}{\}}{0pt}{}{k}{i} \genfrac{[}{]}{0pt}{}{m+1}{k+1} = m! .
\end{equation*}

\item
For $k \geq 0$ and $i \geq 1$, the higher-order Bernoulli numbers $B_k^{(i)}$ are defined by the generating function
\begin{equation*}
\left( \frac{t}{e^t -1}\right)^i = \sum_{k=0}^{\infty} B_k^{(i)} \frac{t^k}{k!},
\end{equation*}
where $B_k^{(1)} = B_k$ are the ordinary Bernoulli numbers. From \cite[Equations (5.2)]{wang} and \eqref{myres2}, we get
\begin{equation}\label{it5}
\sum_{k=0}^{m} (-1)^{k} \genfrac{[}{]}{0pt}{}{m+r+1}{k+r+1}_{r+1} B_k^{(i)} = \binom{m+i}{i}^{-1}
\genfrac{[}{]}{0pt}{}{m+i+r+1}{i+r+1}_{r+1}.
\end{equation}
In particular, for $i=1$ we obtain
\begin{equation*}
\sum_{k=0}^{m} (-1)^{k} \genfrac{[}{]}{0pt}{}{m+r+1}{k+r+1}_{r+1} B_k = \frac{1}{m+1}
\genfrac{[}{]}{0pt}{}{m+r+2}{r+2}_{r+1},
\end{equation*}
which, for $r=0$, reduces to the well-known identity (see, e.g., \cite[Equation (5.6)]{wang})
\begin{equation*}
\sum_{k=0}^{m} (-1)^{k} \genfrac{[}{]}{0pt}{}{m+1}{k+1} B_k = m! H_{m+1}.
\end{equation*}

Moreover, employing the following representation for $B_k^{(i)}$ obtained by Kim {\it et al.\/} \cite{kim}
\begin{equation*}
B_k^{(i)} = \sum_{j=0}^{k} (-1)^j \binom{i+j}{i}^{-1} \genfrac{[}{]}{0pt}{}{i+j}{i} \genfrac{\{}{\}}{0pt}{}{k}{j},
\end{equation*}
and setting $r \to r-1$ in \eqref{it5}, we get
\begin{equation}\label{it6}
\sum_{k=0}^{m} \sum_{j=0}^{k} (-1)^{k+j} \binom{i+j}{i}^{-1} \genfrac{[}{]}{0pt}{}{i+j}{i}
\genfrac{\{}{\}}{0pt}{}{k}{j} \genfrac{[}{]}{0pt}{}{m+r}{k+r}_{r} = \binom{m+i}{i}^{-1}
\genfrac{[}{]}{0pt}{}{m+i+r}{i+r}_{r}.
\end{equation}
In particular, substituting $i=2$ and $r=0$ in \eqref{it6} yields
\begin{equation*}
\sum_{k=0}^{m} \sum_{j=0}^{k} (-1)^{k+j} \frac{j!}{j+2} \genfrac{\{}{\}}{0pt}{}{k}{j}
\genfrac{[}{]}{0pt}{}{m}{k} H_{j+1} = \frac{m!}{m+2} H_{m+1}.
\end{equation*}
On the other hand, employing the representation for $B_k^{(i)}$ given in \cite[Equation (15)]{sriv}
\begin{equation*}
B_k^{(i)} = \sum_{j=0}^{k} (-1)^j \binom{k+i}{i+j} \binom{i+j-1}{i-1} \binom{k+j}{j}^{-1}
\genfrac{\{}{\}}{0pt}{}{k+j}{j},
\end{equation*}
and making $r \to r-1$ in \eqref{it5}, we get
\begin{equation}\label{it7}
\sum_{k=0}^{m} \sum_{j=0}^{k} (-1)^{k+j} \frac{\binom{k+i}{i+j} \binom{i+j-1}{i-1}}{\binom{k+j}{j}}
\genfrac{\{}{\}}{0pt}{}{k+j}{j} \genfrac{[}{]}{0pt}{}{m+r}{k+r}_{r}
= \binom{m+i}{i}^{-1} \genfrac{[}{]}{0pt}{}{m+i+r}{i+r}_{r}.
\end{equation}
In particular, when $i=3$ and $r=0$, \eqref{it7} implies the relation
\begin{equation*}
\sum_{k=0}^{m} \sum_{j=0}^{k} (-1)^{k+j} (j+1)(j+2) \frac{\binom{k+3}{j+3}}{\binom{k+j}{j}}
\genfrac{\{}{\}}{0pt}{}{k+j}{j} \genfrac{[}{]}{0pt}{}{m}{k}
= \frac{6 \cdot m!}{m+3} \big( H_{m+2}^2 - H_{m+2}^{[2]} \big).
\end{equation*}

\item
From \cite[Equations (5.3)]{wang} and \eqref{myres2}, we get
\begin{equation}\label{it8}
\sum_{k=0}^{m} \genfrac{[}{]}{0pt}{}{m+r+1}{k+r+1}_{r+1} B_k(i+1) = m! \binom{m+i+r+1}{i+r}
\big( H_{m+i+r+1} - H_{i+r}\big).
\end{equation}
Taking $r \to r-1$ and $i \to i-1$ in \eqref{it8} gives
\begin{equation*}
\sum_{k=0}^{m} \genfrac{[}{]}{0pt}{}{m+r}{k+r}_{r} B_k(i) = m! \binom{m+i+r-1}{i+r-2}
\big( H_{m+i+r-1} - H_{i+r-2}\big) = m! H_{m+1}^{(i+r-1)},
\end{equation*}
which is just the identity \eqref{fin5} with $x$ replaced by $i$.

\item
In \cite[p. 1508]{wang}, we find the identity
\begin{equation*}
\sum_{k=0}^m P(k,m+r+i,r+i) \frac{B_k^{(i)}}{k!} = \binom{r+i}{i}\binom{m+i}{i}^{-1} P(i,m+r+i,r).
\end{equation*}
Using \eqref{myres2}, we can write this identity in the equivalent form
\begin{equation}\label{1508}
\sum_{k=0}^m  \genfrac{[}{]}{0pt}{}{m+r+i+1}{k+r+i+1}_{r+i+1} B_k^{(i)} = \binom{m+i}{i}^{-1}
\genfrac{[}{]}{0pt}{}{m+i+r+1}{i+r+1}_{r+1}.
\end{equation}
In particular, setting $r=0$ and $i=1$ in \eqref{1508} gives
\begin{equation*}
\sum_{k=0}^m  P_k \big( H(m,1;1), H(m,2;1), \ldots, H(m,k;1) \big) \frac{B_k}{k!} = \frac{H_{m+1}}{m+1}.
\end{equation*}
Incidentally, from \eqref{it5} and \eqref{1508} it follows that
\begin{equation*}
\sum_{k=0}^m  \genfrac{[}{]}{0pt}{}{m+r+i+1}{k+r+i+1}_{r+i+1} B_k^{(i)}
= \sum_{k=0}^{m} (-1)^{k} \genfrac{[}{]}{0pt}{}{m+r+1}{k+r+1}_{r+1} B_k^{(i)}.
\end{equation*}

\item
Combining the recurrence appearing in \cite[p. 1505]{wang}
\begin{equation*}
\binom{j+r}{r} P(i,j+r,r) = \binom{j+r-1}{r-1} P(i,j+r-1,r-1) + \binom{j+r-1}{r} P(i,j+r-1,r)
\end{equation*}
of the numbers $\binom{j+r}{r}P(i,j+r,r)$, and \eqref{myres2}, we readily obtain the corresponding recurrence for the $r$-Stirling numbers of the first kind
\begin{equation*}
\genfrac{[}{]}{0pt}{}{j+r+1}{i+r+1}_{r+1} = \genfrac{[}{]}{0pt}{}{j+r}{i+r}_{r} + j \genfrac{[}{]}{0pt}{}{j+r}{i+r+1}_{r+1}.
\end{equation*}

\end{itemize}

\enlargethispage{2mm}

As a final point to this section, from \cite[Equation (4)]{wuyun} and \eqref{myres2}, we can get the identity
\begin{equation*}
\sum_{j=i}^{m} \frac{1}{m+1-j} \, \frac{1}{(j-1)!} \genfrac{[}{]}{0pt}{}{j+r}{i+r}_{r+1}
= \frac{i}{m!} \genfrac{[}{]}{0pt}{}{m+r+1}{i+r+1}_{r+1},
\end{equation*}
which holds for any integers $1 \leq i \leq m$ and $r \geq 0$. Note the close resemblance of this identity to \eqref{wuyu}. In particular, for $r=0$ we find that
\begin{equation*}
\sum_{j=i}^{m} \frac{1}{m+1-j}\, \frac{1}{(j-1)!} \genfrac{[}{]}{0pt}{}{j}{i}
= \frac{i}{m!} \genfrac{[}{]}{0pt}{}{m+1}{i+1}.
\end{equation*}

\section{Conclusion}

In conclusion, in this note we have shown that the above formulas \eqref{kargin} and \eqref{cere} for the hyper-sums of powers of integers are in fact equivalent, bringing to light an explicit formula for $\genfrac{[}{]}{0pt}{}{m+n+1}{i+n+1}_{n+1}$ as a polynomial in $n$. Moreover, relying on the $r$-Stirling polynomials of the first kind $R_{j+1,i+1}(x)$ and $\overline{R}_{j+1,i+1}(x)$, we have expressed the $i$th derivative of the hyperharmonic polynomials $H_{j+1}^{(x)}$ and $H_{j+1}^{(x+1)}$ as a polynomial in $x$, complementing the formula given in \cite[Equation (19)]{kargin}. Furthermore, we have shown the relationship between the (exponential) complete Bell polynomials and the $r$-Stirling numbers of the first kind, and we have derived some identities involving the Bernoulli numbers and polynomials, the $r$-Stirling numbers of the first kind, the Stirling numbers of both kinds, and the harmonic numbers.

\end{document}